  % This is a Plain TeX file for the paper
  %
  % KMS STATES FOR GENERALIZED GAUGE ACTIONS ON CUNTZ-KRIEGER ALGEBRAS
  % (An application of the Ruelle-Perron-Frobenius Theorem)
  %
  % By Ruy Exel
  %
  % 17-Oct-2001

  % FONTS

  \newcount\fontset
  \fontset=1
  \def\dualfont#1#2#3{\font#1=\ifnum\fontset=1 #2\else#3\fi}

  \dualfont\bbfive{bbm5}{cmbx5}
  \dualfont\bbseven{bbm7}{cmbx7}
  \dualfont\bbten{bbm10}{cmbx10}
  \font \eightbf = cmbx8
  \font \eighti = cmmi8 \skewchar \eighti = '177
  \font \eightit = cmti8
  \font \eightrm = cmr8
  \font \eightsl = cmsl8
  \font \eightsy = cmsy8 \skewchar \eightsy = '60
  \font \eighttt = cmtt8 \hyphenchar\eighttt = -1
  \font \msbm = msbm10
  \font \sixbf = cmbx6
  \font \sixi = cmmi6 \skewchar \sixi = '177
  \font \sixrm = cmr6
  \font \sixsy = cmsy6 \skewchar \sixsy = '60
  \font \tensc = cmcsc10
  
  \font \rs=rsfs10
  \font \titfnt=cmbx12
  \scriptfont \bffam = \bbseven
  \scriptscriptfont \bffam = \bbfive
  \textfont \bffam = \bbten

  \newskip \ttglue

  \def \eightpoint {\def \rm {\fam0 \eightrm }%
  \textfont0 = \eightrm
  \scriptfont0 = \sixrm \scriptscriptfont0 = \fiverm
  \textfont1 = \eighti
  \scriptfont1 = \sixi \scriptscriptfont1 = \fivei
  \textfont2 = \eightsy
  \scriptfont2 = \sixsy \scriptscriptfont2 = \fivesy
  \textfont3 = \tenex
  \scriptfont3 = \tenex \scriptscriptfont3 = \tenex
  \def \it {\fam \itfam \eightit }%
  \textfont \itfam = \eightit
  \def \sl {\fam \slfam \eightsl }%
  \textfont \slfam = \eightsl
  \def \bf {\fam \bffam \eightbf }%
  \textfont \bffam = \eightbf
  \scriptfont \bffam = \sixbf
  \scriptscriptfont \bffam = \fivebf
  \def \tt {\fam \ttfam \eighttt }%
  \textfont \ttfam = \eighttt
  \tt \ttglue = .5em plus.25em minus.15em
  \normalbaselineskip = 9pt
  \def \MF {{\manual opqr}\-{\manual stuq}}%
  \let \sc = \sixrm
  \let \big = \eightbig
  \setbox \strutbox = \hbox {\vrule height7pt depth2pt width0pt}%
  \normalbaselines \rm }

  % HEADER

  \def \Headlines #1#2{\nopagenumbers
    \voffset = 2\baselineskip
    \advance \vsize by -\voffset
    \headline {\ifnum \pageno = 1 \hfil
    \else \ifodd \pageno \tensc \hfil \lcase {#1} \hfil \folio
    \else \tensc \folio \hfil \lcase {#2} \hfil
    \fi \fi }}

  % CONTROL SEQUENCES

  \def \lcase #1{\edef \auxvar {\lowercase {#1}}\auxvar }

  \def \goodbreak {\vskip0pt plus.1\vsize \penalty -250 \vskip0pt
plus-.1\vsize }

  \newcount \secno \secno = 0
  \newcount \stno

  \def \seqnumbering {\global \advance \stno by 1
    \number \secno .\number \stno }

  \def \label #1{\def\localvariable {\number \secno
    \ifnum \number \stno = 0\else .\number \stno \fi }\global \edef
    #1{\localvariable }}

  \def\section #1{\global\def\SectionName{#1}\stno = 0 \global
\advance \secno by 1 \bigskip \bigskip \goodbreak \noindent {\bf
\number \secno .\enspace #1.}\medskip \noindent \ignorespaces}

  \long \def \sysstate #1#2#3{\medbreak \noindent {\bf \seqnumbering
.\enspace #1.\enspace }{#2#3\vskip 0pt}\medbreak }
  \def \state #1 #2\par {\sysstate {#1}{\sl }{#2}}
  \def \definition #1\par {\sysstate {Definition}{\rm }{#1}}

  % Examples
  % \sysstate ...{Theorem}.{font}{Text}
  % \state .......Theorem.........Text\par
  % \definition ..................Text\par

  \def \proof {\medbreak \noindent {\it Proof.\enspace }}
  \def \proofend {\ifmmode \eqno \square \else \hfill \square
\looseness = -1 \medbreak \fi }

  \def\iItem {\smallskip}
  \def\Item #1{\smallskip \item {#1}}
  \newcount \zitemno \zitemno = 0
  \def\izitem {\zitemno = 0}
  \def\zitem {\global \advance \zitemno by 1 \Item {{\rm(\romannumeral
\zitemno)}}}

  \newcount \footno \footno = 1
  \newcount \halffootno \footno = 1
  \def\footcntr {\global \advance \footno by 1
  \halffootno =\footno
  \divide \halffootno by 2
  $^{\number\halffootno}$}
  \def\fn#1{\footnote{\footcntr}{\eightpoint#1}}

  % STANDARD DEFINITIONS

  \def \({\left (}
  \def \){\right )}
  \def \[{\left \Vert }
  \def \]{\right \Vert }
  \def \*{\otimes }
  \def \+{\oplus }
  \def \:{\colon }
  \def \<{\left \langle }
  \def \>{\right \rangle }
  \def \text #1{\hbox {\rm #1}}
  \def \curly#1{\hbox{\rs #1\/}}
  \def \ds{\displaystyle}
  \def \and {\hbox {,\quad and \quad }}
  
  \def \calcat #1{\,{\vrule height8pt depth4pt}_{\,#1}}

  \def \crossproduct {{\hbox {\msbm o}}}
  
  \def \for #1{,\quad \forall\,#1}
  \def \inv {^{-1}}
  \def \pmatrix #1{\left [\matrix {#1}\right ]}
  \def \square {\hbox {$\sqcap \!\!\!\!\sqcup $}}
  \def \stress #1{{\it #1}\/}

  \def \|{\Vert }
  \def \inv {^{-1}}

  % REFERENCE CONTROL

  \newcount \bibno \bibno =0
  \def \newbib #1{\global \advance \bibno by 1 \edef #1{\number
    \bibno}}
  \def\cite #1{{\rm [\bf #1\rm ]}}
  \def\scite #1#2{\cite{#1{\rm \hskip 0.7pt:\hskip 2pt #2}}}
  \def\lcite #1{(#1)}
  \def\fcite #1#2{#1}
  \def\bibitem#1#2#3#4{\smallskip \item {[#1]} #2, ``#3'', #4.}

  \def \references {
    \begingroup
    \bigskip \bigskip \goodbreak
    \eightpoint
    \centerline {\tensc References}
    \nobreak \medskip \frenchspacing }

  % MISCELLANEOUS DEFINITIONS

  \def\N{{\bf N}}
  \def\R{{\bf R}}
  \def\Sp{\Sigma_A}
  \def\Cs{C(\Sp)}
  \def\Tr{{\curly L}}
  \def\ATr{{\cal L}}
  \def\exp{e^}
  \def\compos{\mathop{\raise 1pt \hbox{$\scriptscriptstyle \circ$}}}
  \def\a{\alpha}
  \def\b{\beta}
  \def\d{\delta}
  \def\l{\lambda}
  \def\cp{\mathop{\crossproduct_{\a,\ATr}} \N}
  \def\scp{\mathop{\crossproduct} \N}
  \def\OA{{\cal O}_A}
  \def\On{{\cal O}_n}
  \def\half{^{1/2}}
  \def\Ham{\hbox{\curly H}} \def\Ham{H}
  \def\Hc{H\"older continuous}
  \def\ind{\hbox{ind}}

  % REFERENCES

  \newbib \Baladi
  \newbib \Bowen
  \newbib \BJO
  \newbib \Evans
  \newbib \Amena
  \newbib \Endo
  \newbib \Tower
  \newbib \Kms
  \newbib \OP
  \newbib \QR
  \newbib \RuelleI
  \newbib \RuelleII

  \Headlines {KMS states on Cuntz-Krieger algebras} {Ruy Exel}

  \centerline{\titfnt KMS STATES FOR GENERALIZED GAUGE}
  \medskip
  \centerline{\titfnt ACTIONS ON CUNTZ-KRIEGER ALGEBRAS}
  \medskip
  \centerline{\it (An application of the Ruelle-Perron-Frobenius
Theorem)}

  \bigskip \bigskip
  \centerline
  {\tensc Ruy Exel\footnote{*}{\eightrm Partially supported by CNPq.}}

  \bigskip \bigskip
  \midinsert\narrower\narrower\noindent
  Given a zero-one matrix $A$ we consider certain
one-parameter groups of automorphisms of the Cuntz-Krieger algebra
$\OA$, generalizing the usual gauge group, and depending on a positive
continuous function $\Ham$ defined on the Markov space $\Sp$.  The
main result consists of an application of Ruelle's Perron-Frobenius
Theorem to show that these automorphism groups admit a single KMS
state.
  \endinsert

  \section{Introduction}
  In 1978 Olesen and Pedersen \cite{\OP} showed that the periodic
gauge action on the Cuntz algebra $\On$ admits a unique KMS state,
whose inverse temperature is $\beta = \log n$.  Two years later Evans
\scite{\Evans}{2.2} extended their result to include, among other
things, non-periodic gauge actions, namely one-parameter automorphism
groups on $\On$ given on the standard generating partial isometries
$S_j$ by
  $$
  \gamma_t(S_j) = N_j^{it} S_j
  \for t\in\R,
  $$
  where $\{N_j\}_{j=1}^n$
  is a collection of real numbers with $N_j>1$ for all $j$.  See also
\scite{\BJO}{3.1}.  In 1984 Enomoto, Fujii and Watatani treated the
case of the periodic gauge action on the Cuntz-Krieger algebra $\OA$
for an irreducible matrix $A$ and again arrived at the conclusion that
there exists a unique KMS state.  The case of a non-periodic gauge
action on $\OA$ was discussed in \cite{\Kms} in the context of
Cuntz-Krieger algebras for infinite matrices but, specializing the
conclusions to the finite case, one gets the expected result that if
the matrix $A$ is irreducible and the parameters $N_j$ are all greater
than 1 then there exists a unique KMS state.

The present work aims to take a new step in the direction of
understanding the KMS states on Cuntz-Krieger algebras (over finite
matrices) by studying generalized gauge actions on $\OA$.  In order to
describe these actions let $\Sp$ be the one-sided Markov space for the
given matrix $A$ and consider the copy of $\Cs$ within $\OA$ that is
generated by the elements of the form
  $S_{i_1}\ldots S_{i_k} S_{i_k}^* \ldots S_{i_1}^*$, where the $S_i$
are the standard generating partial isometries.  Fixing an invertible
element $U\in\Cs$ it is not hard to see that the correspondence
  $$
  S_j \mapsto US_j
  $$
  extends to give an automorphism of $\OA$.  Therefore if $H\in\Cs$ is
a strictly positive element there exists a unique one-parameter
automorphism group $\{\gamma_t\}_{t\in\R}$ of $\OA$ such that
  $$
  \gamma_t(S_j) = \Ham^{it}S_j.
  $$
  We will refer to $\gamma$ as the \stress{generalized gauge action}.
  It is easy to see that this in fact generalizes both the periodic
and the non-periodic gauge actions referred to above.

The goal of this paper, as the title suggests, is to study the KMS
states for the generalized gauge action on $\OA$.  Our main result,
Theorem \fcite{4.4}{\MainResult}, states that under certain hypotheses
there exists a single such KMS state.

  The method employed consists of considering $\OA$ as the crossed
product of $\Cs$ by the endomorphism induced by the Markov subshift
\cite{\Endo} and applying Theorem 9.6 from \cite{\Tower} to reduce the
problem to the search for probability measures on $\Sp$ which are
fixed by Ruelle's transfer operator \cite{\RuelleI, \RuelleII,
\Baladi, \Bowen}.  This turns out to be closely related to Ruelle's
version of the Perron-Frobenius Theorem (see
e.g.~\scite{\Bowen}{1.7}), except that the latter deals with
eigenvalues for the transfer operator while we need actual fixed
points.  With not too much effort we are then able to exploit Ruelle's
Theorem in order to understand the required fixed points and thus
reach our conclusion.

It should be stressed that Ruelle's Theorem requires two crucial
hypotheses, namely that the matrix $A$ be \stress{irreducible and
aperiodic} in the sense that there exists a positive integer $m$ such
that all entries of $A^m$ are strictly positive (see
e.g.~\scite{\Baladi}{Section 1.2}), and that $\Ham$ is {\Hc}.  We are
therefore forced to postulate these conditions leaving open the
question as to whether one could do without them.

The organization of this paper is as follows: in section (2), the
longer and more technical section of this work, we give a brief
account of Ruelle's Theorem and draw the conclusions we need with
respect to the existence and uniqueness of probability measures that
are fixed under the transfer operator.

Section (3) is devoted to reviewing results about crossed products by
endomorphisms and in the final section we put all the pieces together
proving our main result.

  I would finally like to acknowledge helpful conversations with
M.~Viana who, among other things, brought Ruelle's Theorem to my
attention.

  \section{Ruelle's Perron-Frobenius Theorem}
  Beyond establishing our notation this section is intended to present
Ruelle's Perron-Frobenius Theorem and to develop some further
consequences of it to be used in later sections.

  Fix, once and for all, an $n\times n$ matrix $A = \{A_{i,j}\}_{1
\leq i,j\leq n}$, with $A_{i,j}\in \{0,1\}$ for all $i$ and $j$, such
that no row or column of $A$ is identically zero.

Throughout this paper we will be concerned with the associated
(one-sided) subshift of finite type, namely the dynamical system
$(\sigma,\Sp)$, where $\Sp$ is the compact topological subspace of the
infinite product space
  $\prod_{i\in \N}\{1,2,\ldots,n\}$ given by
  $$
  \Sp = \Big\{ x = (x_0, x_1, x_2, \ldots\ ) \in \prod_{i\in
\N}\{1,2,\ldots,n\} : A_{x_i,x_{i+1}} = 1 \hbox{ for all } i\ge
0\Big\},
  $$
  and $\sigma : \Sp \to \Sp$ is the ``left shift'', namely the
continuous function given by
  $$
  \sigma(x_0, x_1, x_2, \ldots\ ) = (x_1, x_2, x_3, \ldots\ ).
  $$
  From the assumption that no column of $A$ is identically zero it
follows that $\sigma$ is surjective.

  Given a real number $\beta\in(0,1)$ define a metric $d$ on $\Sp$ by
setting
  $$
  d(x,y) = \beta^{N(x,y)}
  \for x,y\in\Sp,
  $$
  where $N(x,y)$ is the largest integer $N$ such that $x_i=y_i$ for
all $i < N$.  In the special case in which $x=y$ we set
$N(x,y)=+\infty$ and interpret $\beta^{N(x,y)}$ as being zero.  It is
easy to see that this metric is compatible with the product topology.

Let $\Cs$ denote the C*-algebra of all continuous complex functions on
$\Sp$.  We will consider the operator
  $$
  \Tr : \Cs \to \Cs
  $$
  given by
  $$
  \Tr(f)\calcat x =
  \sum_{y\in \sigma\inv(\{x\})} f(y)
  \for f\in \Cs \for x\in \Sp.
  \eqno{(\seqnumbering)}
  \label \DefTransfer
  $$
  Since $\sigma$ is surjective one has that $\sigma\inv(\{x\})$ is
never empty.  It is also clear that $\sigma\inv(\{x\})$ has at most
$n$ elements so that the above sum is finite for every $x$.  One
checks that $\Tr(f)$ is indeed a continuous function and hence that
$\Tr$ is a well defined linear operator on $\Cs$, which is moreover
positive and bounded.

Given a real continuous function $\phi$ on $\Sp$ the operator
  $$
  \Tr_\phi : \Cs \to \Cs
  $$
  given by $\Tr_\phi(f) = \Tr(\exp{\phi}f)$
  was introduced by Ruelle in \scite{\RuelleI}{2.3} (see also
  \cite{\RuelleII}, \cite{\Bowen}, and \cite{\Baladi})
  and it is usually referred to as \stress{Ruelle's transfer
operator}.

Most of the time we will assume that $\phi$ is {\Hc} with respect to
the metric $d$ above: recall that a complex function $\phi$ on a
metric space $M$ is said to be \stress{\Hc} when one can find positive
constants $K$ and $\alpha$ such that
  $
  |\phi(x)- \phi(x)| \leq K d(x,y)^\alpha,
  $
  for all $x$ and $y$ in $M$.

The most important technical tool to be used in this work is the
celebrated Ruelle-Perron-Frobenius Theorem which we now state for the
convenience of the reader.

  \state Theorem
  \label \RuelleThm
  {\rm (D.~Ruelle)} Let $A$ be an $n\times n$ zero-one matrix and let
$\phi$ be a real function defined on $\Sp$.  Suppose that:
  \iItem
  \Item{(a)} There exists a positive integer $m$ such that $A^m>0$ (in
the sense that all entries are $>0$), and
  \Item{(b)} $\phi$ is {\Hc}.
  \medskip \noindent Then there are:
  a strictly positive function $h\in\Cs$,
  a Borel probability measure $\nu$ on $\Sp$, and
  a real number $\l>0$,
  such that
  \izitem
  \zitem $\Tr_\phi(h) = \l h$,
  \zitem $\Tr_\phi^*(\nu) = \l \nu$, where $\Tr_\phi^*$ is the adjoint
operator acting on the dual of\/ $\Cs$, and
  \zitem for every $g\in\Cs$ one has that
  $
  \ds\lim_{k\to\infty} \|\l^{-k} \Tr_\phi^k(g) - \nu(g)h\| = 0.
  $
  \proof See e.g.~\scite{\Bowen}{1.7}.

  \state Proposition
  \label \UniqueLambdaNu
  Under the hypotheses of \lcite{\RuelleThm} there exists a unique
pair $(\l_1,\nu_1)$ such that
  $\l_1$ is a complex number,
  $\nu_1$ is a probability measure on $\Sp$,
  and
  $\Tr_\phi^*(\nu_1) = \l_1 \nu_1$.

  \proof The existence obviously follows from \lcite{\RuelleThm.ii}.
As for uniqueness let $(\l_1,\nu_1)$ be such a pair and let $(\l,\nu)$
be as in \lcite{\RuelleThm}.  For all $g\in\Cs$ we have
  $$
  \lim_{k\to\infty} \({\l_1 \over \l}\)^k\nu_1(g) =
  \lim_{k\to\infty} \nu_1\(\l^{-k}\Tr_\phi^k(g)\) =
  \nu(g)\nu_1(h),
  $$
  by \lcite{\RuelleThm.iii}.  Plugging $g=1$ above we conclude that
the sequence $\({\l_1 \over \l}\)^k$ converges to the nonzero value
$\nu_1(h)$ but this is only possible if $\l_1=\l$.  For every $g$ we
then have that
  $
  \nu_1(g) =
  \nu(g)\nu_1(h),
  $
  so $\nu_1$ is proportional to $\nu$.  But since these are
probability measures we must have $\nu_1=\nu$.  \proofend

  In particular it follows that both the $\l$ and the $\nu$ in the
conclusion of \lcite{\RuelleThm} are uniquely determined.  In the
following we give an explicit way to compute $\l$ in terms of
$\Tr_\phi$ (see \scite{\Baladi}{1.39}).

  \state Proposition
  \label \BetaLambda
  Under the hypotheses of \lcite{\RuelleThm} one has that
  $$
  \l = \lim_{k\to\infty} \|\Tr^k_\phi(1)\|^{1/k}.
  $$

  \proof
  Plugging $g=1$ in \lcite{\RuelleThm.iii} we conclude that
  $\ds
  \lim_{k\to\infty} \l^{-k}\|\Tr_\phi^k(1)\| =
  \|h\| >0.
  $
  So we may choose $n_0\in\N$ such that for all $n\geq n_0$
  $$
  {\|h\|\over2} <
  \l^{-k}\|\Tr_\phi^k(1)\| <
  2\|h\|.
  $$
  Taking $k^{th}$ roots and then the limit as $k\to\infty$ we get the
conclusion.
  \proofend

In the application of Ruelle's Theorem that we have in mind we will
take
  $$
  \phi = \phi_\b = -\b\log(\Ham),
  \eqno{(\seqnumbering)}
  \label \PhiBeta
  $$
  where $\Ham$ is a strictly positive continuous function on $\Sp$ and
$\b>0$ is a real number.

Observe that if $\Ham$ is {\Hc} then so is $\phi_\b$ for every real
$\b$ (this is because ``$\log$'' is Lipschitz on every compact subset
of $(0,+\infty)$, e.g.~the range of $\Ham$).  In this case Ruelle's
Theorem gives a correspondence $\b\to\l$ which we would like to
explore more closely in what follows.

  \state Proposition
  \label \MyContribution
  Let
  $A$ be an $n\times n$ zero-one matrix satisfying
\lcite{\RuelleThm.a}
  and suppose that
  $\Ham$ is a {\Hc} function on $\Sp$ such that
  $$
  \Ham(y) > 1
  \for y\in\Sp.
  $$
  For every $\b\geq0$
  let $\phi_\b$ be as in \lcite{\PhiBeta} and denote by $\l(\b)$ the
unique $\l$ satisfying the conditions of \lcite{\RuelleThm} for
$\phi=\phi_\b$.  Then one has that
  \izitem
  \zitem $\l(0) > 1$,
  \zitem $\ds\lim_{\b\to\infty}\l(\b) =0$, and
  \zitem $\l$ is a strictly decreasing continuous function of $\b$.

  \proof
  Observe that
  $$
  \Tr_{\phi_\b}(f) = \Tr(\exp{\phi_\b}f) =
  \Tr\(\exp{-\b\log(\Ham)}f\) =
  \Tr\(\Ham^{-\b}f\).
  $$
  Let $m$ and $M$ be the supremum and infimum of $\Ham$ on $\Sp$,
respectively.
  For every $\b\geq 0$ and $y\in \Sp$ one therefore has that
  $$
  M^{-\b} \leq \Ham(y)^{-\b} \leq m^{-\b},
  $$
  so that if $f\in\Cs$ is nonnegative we have
  $$
  M^{-\b}\Tr(f) \leq \Tr_{\phi_\b}(f) \leq m^{-\b} \Tr(f).
  $$
  By induction it is easy to see that for all $k\in \N$
  $$
  M^{-k\b}\Tr^k(f) \leq \Tr_{\phi_\b}^k(f) \leq m^{-k\b} \Tr^k(f).
  $$
  Taking norms and $k^{th}$ roots we conclude that
  $$
  M^{-\b}\|\Tr^k(f)\|^{1/k} \leq \|\Tr_{\phi_\b}^k(f)\|^{1/k} \leq
  m^{-\b} \|\Tr^k(f)\|^{1/k}.
  $$
  Plugging $f=1$ above and observing that
  $1 \leq \|\Tr^k(1)\| \leq \|\Tr\|^k$
  we obtain
  $$
  M^{-\b} \leq \|\Tr_{\phi_\b}^k(1)\|^{1/k} \leq m^{-\b} \|\Tr\|,
  $$
  and hence \lcite{\BetaLambda} yields
  $$
  M^{-\b} \leq \l(\b) \leq m^{-\b} \|\Tr\|.
  $$
  Observing that $\Ham>1$, and hence that $m>1$, we deduce (ii).  It
is also clear from the above that $\l(0) \geq 1$ so it is enough to
show that $\l(0) \neq 1$ in order to obtain (i).

Arguing by contradiction suppose that $\l(0)=1$.  Let $h>0$ be given
by \lcite{\RuelleThm} so that $\Tr_{\phi_0}(h) = \Tr(h) = h$.  Choose
$x^0\in\Sp$ such that
  $
  h(x^0) = \inf_{y\in \Sp} h(y),
  $
  and observe that, since
  $$
  h(x^0) = \sum_{y\in \sigma\inv(\{x^0\})} h(y),
  $$
  there exists a unique $y$ in $\sigma\inv(\{x^0\})$ which moreover
satisfies $h(y) = h(x^0)$.  Repeating this process one obtains a
sequence $\{x^k\}_{k\in\N}$ in $\Sp$ such that
  $\sigma\inv(\{x^k\}) = \{x^{k+1}\}$ for all $k$.  Letting $x_k =
x^k_0$ (the zero$^{th}$ coordinate of $x^k$) we have that
  $A_{x_{k+1}, x_k}=1$ and also that this is the only nonzero entry of
$A$ in the column $x_k$.  Since $A$ is a finite matrix the sequence
$\{x_k\}$ must be periodic.  Assuming without loss of generality that
the first period of this sequence is $\{1,\ldots,m\}$, where $m\leq
n$, we see that $A$ has the form
  $$
  A = \pmatrix{ S_m & B \cr 0 & C},
  $$
  where $S_m$ is the matrix of the forward permutation of $m$
elements.  However this is easily seen to contradict
\lcite{\RuelleThm.a} both when $m<n$ (because the zero block in the
lower left corner will appear in any power of $A$) and when $m=n$
(because $S_m$ definitely fails to satisfy \lcite{\RuelleThm.a}).

  In order to prove (iii)
  let $\d>0$ so that $m^\d \leq \Ham(y)^\d \leq M^\d$ for all $y$ in
$\Sp$.  Given $\b\in\R$ we then have that
  $$
  m^\d \Ham(y)^{-\b} \leq
  \Ham(y)^{-(\b-\d)} \leq
  M^\d \Ham(y)^{-\b}.
  $$
  For every nonnegative continuous function $f$ it follows that
  $$
  m^\d\Tr_{\phi_\b}(f) \leq
  \Tr_{\phi_{\b-\d}}(f) \leq
  M^\d\Tr_{\phi_\b}(f),
  $$
  and if $k\in\N$ one has
  $$
  m^{k\d}\Tr_{\phi_\b}^k(f) \leq
  \Tr_{\phi_{\b-\d}}^k(f) \leq
  M^{k\d}\Tr_{\phi_\b}^k(f).
  $$
  Taking norms and $k^{th}$ roots we conclude that
  $$
  m^\d\|\Tr_{\phi_\b}^k(f)\|^{1/k} \leq
  \|\Tr_{\phi_{\b-\d}}^k(f)\|^{1/k} \leq
  M^\d\|\Tr_{\phi_\b}^k(f)\|^{1/k}.
  $$
  With $f=1$ and taking the limit as $k\to\infty$, we get by
\lcite{\BetaLambda} that
  $$
  m^\d\l(\b) \leq
  \l(\b-\d) \leq
  M^\d\l(\b).
  \eqno{(\seqnumbering)}
  \label \IneqOne
  $$
  Substituting $\b+\d$ for $\b$ above leads to
  $$
  M^{-\d}\l(\b) \leq
  \l(\b+\d) \leq
  m^{-\d}\l(\b).
  \eqno{(\seqnumbering)}
  \label \IneqTwo
  $$
  By \lcite{\IneqOne} and \lcite{\IneqTwo} one sees that $\l$ is a
continuous function of $\b$.  Since $m>1$ by hypothesis the rightmost
inequality in \lcite{\IneqTwo} gives
  $\l(\b+\d) <\l(\b)$ and hence that $\l$ is strictly decreasing.
  \proofend

  \state Corollary
  \label \LambdaIsOne
  Under the hypotheses of \lcite{\MyContribution} there exists a
unique $\b >0$ such that $\lambda(\b)=1$.

  \section{Preliminaries on Crossed Products}
  Define the map
  $\alpha:\Cs \to \Cs$ by the formula
  $$
  \alpha(f) = f \compos \sigma
  \for f\in \Cs.
  $$
  It is easy to see that $\alpha$ is a C*-algebra endomorphism of
$\Cs$.  Since $\sigma$ is surjective one has that $\a$ is injective.
We should also notice that $\a(1)=1$.

  For $x\in\Sp$ let
  $$
  Q(x) = \#\big\{y\in X: \sigma(y)=x\big\},
  $$
  (``$\#$'' meaning number of elements).
  Alternatively $Q(x)$ may be defined as the number of ``ones'' in the
column of $A$ indexed by $x_0$.  Therefore $1 \leq Q(x) \leq n$ for
all $x\in\Sp$ so that in particular $Q$ is invertible as an element of
$\Cs$.
  Define the operator
  $$
  \ATr : \Cs \to \Cs
  $$
  by $\ATr(f) = Q\inv \Tr(f)$, where $\Tr$ is defined in
  \lcite{\DefTransfer}.  It is easy to see that $Q = \Tr(1)$ and hence
that
  $\ATr(1)=1$.  Moreover
  $$
  \ATr\big(\a(f)g\big) = f\ATr(g)\for f,g\in \Cs,
  $$
  which tells us that $\ATr$ is a \stress{transfer operator} for the
pair $(\Cs,\a)$ according to Definition \lcite{2.1} in \cite{\Endo}.
  One may therefore construct the crossed product algebra
  $$
  \Cs\cp,
  $$
  or $\Cs\scp$, for short, as in \scite{\Endo}{3.7}, which turns out
to be a C*-algebra generated by a copy of $\Cs$ and an extra element
$S$ which, among other things, satisfies
  \iItem
  \Item{$\bullet$} $S^*S= 1$,
  \Item{$\bullet$} $S f = \a(f) S$, and
  \Item{$\bullet$} $S^* f S = \ATr(f),$

  \medskip\noindent
  for all $f\in \Cs.$
  See \cite{\Endo} for the precise definition of $\Cs\scp$.

  In \scite{\Endo}{6.2} it is proved that $\Cs\scp$ is isomorphic to
the Cuntz-Krieger algebra $\OA$.
  It will be convenient for us to bear in mind the isomorphism between
$\OA$ and $\Cs\scp$ given in \cite{\Endo}, which we next describe.
For this consider for each $j=1,\ldots,n$, the clopen subset
$\Sigma^j$ of $\Sp$ given by
  $$
  \Sigma^j = \{x\in\Sp : x_0=j\}.
  $$
  These are precisely the sets forming the standard Markov partition
of $\Sp$.  Also let $P_j$ be the characteristic function of
$\Sigma^j$.  According to \cite{\Endo} there exists an isomorphism
  $$
  \Psi : \OA \to \Cs\scp
  $$
  which is determined by the fact that the canonical generating
partial isometries $S_j\in\OA$ are mapped under $\Psi$ as follows:
  $$
  \Psi(S_j) = P_j S Q\half.
  $$

  We would next like to review the definition of the generalized gauge
action on $\OA$.  For this fix a strictly positive element
$\Ham\in\Cs$.  According to \scite{\Tower}{6.2} there exists a unique
one parameter automorphism group
  $\gamma$
  of $\Cs\scp$ such that for all $t\in\R$,
  $$
  \gamma_t(S) = \Ham^{it}S
  \and
  \gamma_t(f) = f
  \for f\in \Cs.
  $$
  Transferring $\gamma$ to $\OA$ via the isomorphism $\Psi$ described
above one gets an automorphism group on $\OA$ which is characterized
by the fact that
  $$
  \gamma_t(S_j) = \Ham^{it}S_j
  \for j=1,\ldots,n.
  $$

  Observe that in case $\Ham$ is a constant function, say everywhere
equal to Neper's number $e$, and $A_{ij}\equiv 1$, then $\OA$
coincides with the Cuntz algebra $\On$ and one recovers the action
over $\On$ considered in \cite{\OP}.  We shall refer to this as the
\stress{scalar gauge action}.

For a slightly more general example suppose that $\Ham$ is constant on
each $\Sigma^j$, taking the value $N_j$ there.  Then
  $$
  \gamma_t(S_j) = \Ham^{it}S_j = \Ham^{it}P_jS_j = N_j^{it}S_j,
  $$
  and we obtain special cases of actions studied in \cite{\Evans} or
\cite{\Kms}.

  Observe that the composition $E= \a\compos\ATr$ is a conditional
expectation from $\Cs$ onto the range of $\a$.  By
\scite{\Tower}{Section 11}, using the set $\{P_1,\ldots,P_n\}$, we see
that $E$ is of index-finite type.  It therefore follows from
\scite{\Tower}{8.9} that there exists a unique conditional expectation
  $$
  G : \OA \to \Cs
  $$
  which is invariant under the scalar gauge action.  This conditional
expectation must therefore coincide with the conditional expectation
given by \scite{\Amena}{2.9} for the Cuntz-Krieger bundle (see
\cite{\QR} and \cite{\Amena}).

Let us now give a concrete description of $G$ based on the well known
fact that $\OA$ is linearly spanned by the set of all $S_\mu S_\nu^*$,
where $\mu$ and $\nu$ are finite words in the alphabet
$\{1,\ldots,n\}$, and we let $S_\mu = S_{\mu_0}\ldots S_{\mu_k}$
whenever $\mu=\mu_0\ldots\mu_k$.

For any such $\mu$ and $\nu$ we have by \cite{\Amena} that
  $$
  G(S_\mu S_\nu^*) =
  \left\{\matrix{S_\mu S_\nu^* &, \hbox{ if }\mu=\nu, \cr
                 0&, \hbox{ if }\mu\neq\nu} \right.
  \eqno{(\seqnumbering)}
  \label \CondExp
  $$

  \section{KMS states}
  It is our main goal to describe the KMS states on $\OA$ for the
gauge action $\gamma$ determined by a given $\Ham$ as above.  Recall
from \scite{\Tower}{9.6} that for every $\b>0$ the correspondence
  $$
  \psi \mapsto \nu = \psi|_{\Cs}
  \eqno{(\seqnumbering)}
  \label \Correspondence
  $$
  is a bijection from
  the set of KMS$_\b$ states $\psi$ on $\Cs\scp$ and
  the set of probability measures\fn{By the Riesz Representation
Theorem we identify probability measures and states as usual.}  $\nu$
on $\Sp$ such that
  $$
  \nu(f) = \nu\Big(\ATr\big(\Ham^{-\b}\ind(E) f\big)\Big)
  \for f\in \Cs,
  \eqno{(\seqnumbering)}
  \label \MainCondition
  $$
  where $\ind(E)$ is the Jones-Kosaki-Watatani index of $E$.  See
\cite{\Tower} for details.  As observed in \scite{\Tower}{Section 11}
the right hand side of \lcite{\MainCondition} coincides with
  $
  \nu\big(\Tr_{\phi_\b}(f)\big),
  $
  where $\phi_\b$ is as in \lcite{\PhiBeta}, so that
\lcite{\MainCondition} is equivalent to
  $$
  \Tr_{\phi_\b}^*(\nu) = \nu.
  \eqno{(\seqnumbering)}
  \label \MainCondition
  $$

  We now arrive at our main result.

  \state Theorem
  \label \MainResult
  Let $A$ be an $n\times n$ zero-one matrix satisfying
\lcite{\RuelleThm.a}
  and let
  $\Ham$ be a {\Hc} function on $\Sp$ such that
  $
  \Ham(y) > 1
  $
  for all $y$ in $\Sp$.
  Let $\gamma$ be the unique one-parameter automorphism group of $\OA$
such that
  $$
  \gamma_t(S_j) = \Ham^{it}S_j
  \for j=1,\ldots,n,
  $$
  where the $S_j$ are the canonical partial isometries generating
$\OA$.  Then $\Cs\scp$ admits a unique KMS state $\psi$ for $\gamma$.
The inverse temperature at which this state occurs is the unique value
of $\b$ for which $\l(\b)=1$ (see \LambdaIsOne).  In addition $\psi$
is given by
  $$
  \psi = \nu\compos G,
  $$
  where $G$ is the conditional expectation described in
\lcite{\CondExp} and $\nu$ is the unique measure on $\Sp$ satisfying
the Ruelle-Perron-Frobenius Theorem for $\phi= -\b\log(\Ham)$.
Finally there are no ground states for $\gamma$.

  \proof
  By \lcite{\LambdaIsOne} let $\b>0$ be such that $\l(\b)=1$.
Applying \lcite{\RuelleThm} for $\phi=\phi_\b = -\b\log(\Ham)$, let
$\nu$ be the unique probability measure on $\Sp$ satisfying
  $$
  \Tr_{\phi_\b}^*(\nu) = \l(\b)\nu = \nu.
  $$ Then \lcite{\MainCondition} holds and hence by
\scite{\Tower}{9.6} the composition $\psi = \nu\compos G$ is a
KMS$_\b$ state for $\gamma$.

Suppose now that $\b_1>0$ and let $\psi_1$ be a KMS$_{\b_1}$ state for
$\gamma$.  Set $\nu_1 = \psi_1|_{\Cs}$ and observe that, again by
\scite{\Tower}{9.6}, one has that $\nu_1$ satisfies
\lcite{\MainCondition} for $\phi_{\b_1}$.  So the pair $(1,\nu_1)$
satisfies the conditions of \lcite{\UniqueLambdaNu} and hence
$\l(\b_1)=1$ so that $\b_1=\b$ by \lcite{\LambdaIsOne}.  Also by
\lcite{\UniqueLambdaNu} $\nu_1$ must coincide with $\nu$ and hence
$\psi_1=\psi$ because the correspondence in \lcite{\Correspondence} is
bijective.

  That no ground states exist follows from \scite{\Tower}{10.1}.
  \proofend

  \endgroup

\bigskip 
\bigskip
\hfill 17 October 2001

\bigskip 
\obeylines 
\tensc
Departamento de Matem\'atica
Universidade Federal de Santa Catarina
88040-900 Florian\'opolis SC
Brazil

\bigskip
\halign{\indent  #&{\rm #}\cr
E-mail: &exel@mtm.ufsc.br\cr
&exel@ime.usp.br\cr}

  \end